\documentclass[11pt]{amsart}

\usepackage[a4paper,hmargin=3.5cm,vmargin=4cm]{geometry}
\usepackage{amsfonts,amssymb,amscd,amstext}
\usepackage{graphicx,float}
\usepackage[dvips]{epsfig}
\usepackage{pstricks,pst-plot}
\usepackage{hyperref}

\usepackage{fancyhdr}
\input xy
\xyoption{all}

\usepackage{enumerate}
\usepackage{titlesec}
\usepackage{mathrsfs}

\titleformat{\section}%[display]
{\filcenter\bfseries\large} {\thesection{.}}{0.2cm}{}%[$\vspace*{-1.0cm}$]
\titleformat{\subsection}[runin]
{\bfseries} {\thesubsection{.}}{0.15cm}{}[.]
\titleformat{\subsubsection}[runin]
{\em}{\thesubsubsection{.}}{0.15cm}{}[.]
\usepackage[up,bf]{caption}
\theoremstyle{definition}

\newcommand\Gcal{\mathcal{G}}

\newcommand\Zscr{\mathscr{Z}}

\def\c{\mathbb{C}}

\def\n{\mathbb{N}}
\renewcommand\r{\mathbb{R}}

\newcommand\igot{\mathfrak{i}}

\renewcommand\igot{\mathfrak{i}}

\renewcommand\imath{\igot}

\newcommand\di{\partial}

\usepackage{color}

\begin{document}

\fancyhead[LO]{Erratum} 
\fancyhead[RE]{A.\ Alarc\'on, I.\ Castro-Infantes,  and F.\ J.\ L\'opez}
\fancyhead[RO,LE]{\thepage}

\thispagestyle{empty}

%% Title
%\vspace*{7mm}
\begin{center}
{\bf \Large Erratum to ``Interpolation and optimal hitting for complete minimal surfaces with finite total curvature"}

\vspace*{5mm}

%% Authors
{%\large
\bf Antonio Alarc\'on, Ildefonso Castro-Infantes, and Francisco J.\ L\'opez}
\end{center}

%\vspace*{7mm}

%\begin{quote}
%{\small
%\noindent {\bf Abstract} \hspace*{0.1cm}
%In this note we correct a mistake in the article \cite{AlarconCastro-InfantesLopez2019CVPDE}, by the authors, with the reference: Calc. Var. (2019) 58: 21, DOI:10.1007/s00526-018-1465-0. We use the notation from \cite{AlarconCastro-InfantesLopez2019CVPDE}.
%\medskip
%
%\noindent{\bf Keywords} \hspace*{0.1cm} 
%minimal surface, finite total curvature, Riemann surface, meromorphic function, interpolation theory, optimal hitting.
%
%\medskip
%
%\noindent{\bf MSC (2010)} \hspace*{0.1cm} 
%53A10, %Minimal surfaces
%52C42, %Immersions (minimal, prescribed curvature, tight, etc.)
%30D30, % Meromorphic functions, general theory
%32E30. %Holomorphic and polynomial approximation, Runge pairs, interpolation
%}
%\end{quote}

%%%%%%%%%%
%%%%%%%%%%
%%%%%%%%%%
%%%%%%%%%%  Text
%%%%%%%%%%
%%%%%%%%%%

Throughout this note we use the notation from \cite{AlarconCastro-InfantesLopez2019CVPDE}. 

Franc Forstneri\v c pointed out that \cite[Proposition 2.3]{AlarconCastro-InfantesLopez2019CVPDE} is incorrect by giving us a counterexample; we are grateful to him for reporting and for helpful discussions. We overlooked the possible contribution of the cusp points of the curve $\gamma$ in the result to the absolute value of its winding number. 
To fix this, we provide a corrected version of \cite[Proposition 2.3]{AlarconCastro-InfantesLopez2019CVPDE} which measures the possible influence of the cusp points on the winding number of the curve (see Proposition 2.3). We then use this result to obtain corrected versions of \cite[Theorems 1.3 and 1.4]{AlarconCastro-InfantesLopez2019CVPDE} on optimal hitting theory for complete minimal surfaces (see Theorems 1.3 and 1.4). The other main results in \cite{AlarconCastro-InfantesLopez2019CVPDE}, namely \cite[Theorem 1.1 and 3.1]{AlarconCastro-InfantesLopez2019CVPDE}, are not affected by the error. 

Here are the corrected statements; see \cite[Sec.\ 1 and Subsec.\ 2.5]{AlarconCastro-InfantesLopez2019CVPDE} for background.

%\begin{proposition}\label{pro:2.3}
\noindent{\bf Proposition 2.3.} {\em Let $\gamma(t)$ be a real analytic, piecewise regular, closed curve in $\r^2$ admitting a regular normal field. Denote by $t_\gamma$ and $m_\gamma$ the turning number and the number of cusp points of $\gamma$, respectively. Fix $p\in\r^2\setminus\gamma$ and let $w_\gamma(p)$ denote the winding number of $\gamma$ with respect to $p$. Then we have $|w_\gamma(p)|\le 2t_\gamma+\frac12m_\gamma$.}
%\end{proposition}

The addend $\frac12m_\gamma$ is missing in the bound stated in \cite[Proposition 2.3]{AlarconCastro-InfantesLopez2019CVPDE}. 

%\begin{theorem}\label{th:1.3}
\noindent{\bf Theorem 1.3.} 
{\em Let $r\ge 1$ be an integer. For any integer $m$ with  $2-2r \le m\le 1$  there is a set $A_{r;m}\subset\r^3$ which is against the family  $\bigcup_{k\le m}\Zscr_{r;k}$ and consists of $12r+50r^3+2m+1$ points whose affine span is a plane. In particular, the set $A_{r;1}$, which consists of $12r+50r^3+3$ points, is against the family $\Zscr_r$. 

Thus, if $X\colon M\to\r^3$ is a complete nonflat orientable immersed minimal surface with empty boundary and the Euler characteristic $\chi(M)\le m$, and  if $A_{r;m}\subset X(M)$, then the total curvature ${\rm TC}(X)<-4\pi r$. In particular, no complete nonflat orientable immersed minimal surface $X$ with $|{\rm TC}(X)|\le 4\pi r$ contains $A_{r;1}$.}
%\end{theorem}

%\begin{theorem}\label{th:1.4}
\noindent{\bf Theorem 1.4.}
{\em Let $X\colon M\to \r^3$ be a complete orientable immersed minimal surface of finite total curvature and empty boundary. If $L\subset\r^3$ is a straight line  not contained in $X(M)$, then
\[
	\# \big(X^{-1}(L)\big)\leq 6 {\rm Deg}(N)+25{\rm Deg}(N)^3+\chi(M)
\]
where ${\rm Deg}(N)$ is the degree of the Gauss map $N$ of $X$ and $\chi(M)$ is the Euler characteristic of $M$.
}
%\end{theorem}

The terms $50r^3$ and $25 {\rm Deg}(N)^3$ are not present in the statements of \cite[Theorems 1.3 and 1.4]{AlarconCastro-InfantesLopez2019CVPDE}, respectively. They come from the term $\frac12m_\gamma$ in Proposition 2.3.

%\smallskip
\noindent{\em Proof of Proposition 2.3.} 
We follow \cite[proof of Proposition 2.3]{AlarconCastro-InfantesLopez2019CVPDE} until the very end. Note that each curve $\gamma_j$ there, $j=1,\ldots,2t_\gamma$, has precisely $m_j-1$ cusp points, and hence $m_\gamma=\sum_{j=1}^{2t_\gamma}(m_j-1)$. Since every two consecutive subarcs $\gamma_{j,i}$, $\gamma_{j,i+1}$ of $\gamma_j$   have the opposite character (i.e., one is positive and the other one negative) and $\gamma_j$ consists of precisely $m_j$ arcs $\gamma_{j,i}$, it turns out that the absolute value of the signed number of crossings of $\gamma_j$ with $\ell$ is at most $1+E(\frac{m_j-1}2)$, where $E(\cdot)$ denotes integer part, and hence $|w_\gamma(p)|\leq\sum_{j=1}^{2 t_\gamma} \big(1+E(\frac{m_j-1}2)\big)\leq 2 t_\gamma+\frac12 m_\gamma$.
\qed
%\smallskip

\noindent
{\em Proof of Theorem 1.4.}
We assume without loss of generality that $X$ is primitive in the sense that there is no nontrivial (finite) holomorphic covering $\rho\colon M\to M_0$ into an open Riemann surface $M_0$ admitting a conformal minimal immersion $X_0\colon M_0\to\r^3$ such that $X=X_0\circ \rho$.
Also, we assume that $X$ is neither a plane nor a catenoid; otherwise the conclusion of the theorem is obvious.  
The proof then follows word by word that of \cite[Theorem 1.4]{AlarconCastro-InfantesLopez2019CVPDE} until the top of page 18 where  \cite[Proposition 2.3]{AlarconCastro-InfantesLopez2019CVPDE} is applied. Instead, we apply Proposition 2.3 and obtain
$\Big|\sum_{j=1}^m f_*(c_j)\Big|\le \sum_{j=1}^m |w_j| \le  \sum_{j=1}^m (2t_j+\frac12m_j)=2 \sum_{j=1}^m t_j+\frac12\sum_{j=1}^m m_j$, where $m_j$ denotes the number of cusp points of the curve $\alpha_j$ (this replaces (4.5) in \cite{AlarconCastro-InfantesLopez2019CVPDE}). As in the proof in \cite{AlarconCastro-InfantesLopez2019CVPDE}, we have $\sum_{j=1}^m t_j={\rm Deg}(N|_{\overline\Omega})$, and so if we set $m_\Omega=\sum_{j=1}^m m_j$
then
$\Big|\sum_{j=1}^m f_*(c_j)\Big| \le  2{\rm Deg}(N|_{\overline\Omega})+\frac12m_\Omega
$ (this replaces (4.6) in \cite{AlarconCastro-InfantesLopez2019CVPDE}). Combining %\eqref{eq:4.6} 
this with equations (4.2), (4.3), and (4.4) in \cite{AlarconCastro-InfantesLopez2019CVPDE}, we get that
$\#(X^{-1}(L)\cap\overline\Omega)=\#(X^{-1}(L)\cap\Omega)\le2{\rm Deg}(N|_{\overline \Omega})+\frac12m_\Omega+ \sum_{j=1}^r I_{q_j}$.
Denote by $\Omega_1,\ldots,\Omega_a$ $(a\in\n)$ the connected components of $\Sigma\setminus Q$. Joining together the above information for all the $\Omega_i$'s, and taking into account that each Jordan curve in $Q$ lies in the boundary of exactly two of them, we obtain that
$\#(X^{-1}(L))\le 4{\rm Deg}(N)+\frac12\sum_{i=1}^am_{\Omega_i}+\sum_{q\in E}I_q$
(this replaces (4.7) in \cite{AlarconCastro-InfantesLopez2019CVPDE}). The Jorge-Meeks formula (see (2.5) in \cite{AlarconCastro-InfantesLopez2019CVPDE}) then gives 
\begin{equation}\label{eq:mOmega}
	\#(X^{-1}(L))\le 6{\rm Deg}(N)+\chi(M)+\frac12\sum_{i=1}^am_{\Omega_i}.
\end{equation}
To complete the proof we shall now provide an upper bound for $\sum_{i=1}^am_{\Omega_i}$ in terms of the degree of the Gauss map $N$ of $X$. This term is precisely what is missing in \cite[proof of Theorem 1.4]{AlarconCastro-InfantesLopez2019CVPDE} and what makes the addend $25{\rm Deg}(N)^3$ to appear in the statement of Theorem 1.4. Indeed, we claim that
 \begin{equation}\label{eq:cuspglobal-FTC}
	 \sum_{i=1}^a m_{\Omega_i}\leq   50 {\rm Deg}(N)^3.
\end{equation}
To prove this, denote by $g$ the complex Gauss map of $X=(X_1,X_2,X_3)$.  
Since the affine line $L\subset\r^3$ is assumed to be the $x_3$-axis and to lie in $\Gcal_0$ (see conditions {\rm (a)}, {\rm (b)}, and {\rm (c)} in \cite[p.\ 15]{AlarconCastro-InfantesLopez2019CVPDE}), we have that $Q=Q_L=|g|^{-1}(1)=\bigcup_{i=1}^a  b\Omega_i$. Denote by $C$ the set of points in $|g|^{-1}(1)$ whose image by the map $(X_1,X_2)\colon |g|^{-1}(1)\to\r^2$ is a singular point. 
Note that $C$ is finite and contains all cusp points of $(X_1,X_2)|_{|g|^{-1}(1)}$, and hence to prove \eqref{eq:cuspglobal-FTC} it suffices to show that 
 \begin{equation}\label{eq:cuspglobal-FTC-C}
 	\# C\le 25{\rm Deg}(g)^3;
\end{equation} 
note that ${\rm Deg}(g)={\rm Deg}(N)$ and that each Jordan curve in $|g|^{-1}(1)$ lies in the boundary of exactly two domains $\Omega_i$. 
For, observe that the meromorphic function $F=g\frac{\di X_3}{dg}$
on $\Sigma=M\cup E$ satisfies
\begin{equation}\label{eq:C-m_Omega}
	C=\{p\in |g|^{-1}(1)\colon \Re F(p)=0\}.
\end{equation}
Choose $\theta\in \c$ such that $|\theta|=1$  and $g^{-1}(\theta)$ is disjoint from $C$, consider the meromorphic function $u=\imath \frac{g+\theta}{g -\theta}$ on $\Sigma$, and notice that 
 \begin{equation}\label{eq:C-m_Omega1}
	 |g|^{-1}(1)=u^{-1}(\r\cup\{\infty\})\quad \text{and}\quad C\subset u^{-1}(\r).
\end{equation}
Since $X$ is not a plane or a catenoid, we have that $u$ and $F$ are nonconstant, and hence \cite[Proposition IV.11.6]{FarkasKra1992} provides a nonconstant irreducible complex polynomial  $P$ in two variables such that $P(u,F)=0$ everywhere on $\Sigma$. 
Since the Weierstrass data of $X\colon M=\Sigma\setminus E\to\r^3$ are determined by $g$, $F$, and $d g$ and the immersion $X$ is primitive, it is not difficult to see that $(u,F)$ is a primitive pair on $\Sigma$ in the sense of \cite[p.\ 248]{FarkasKra1992}, and hence $\Sigma$ is identified with the algebraic curve associated to $\{P=0\}$ via the biholomorphism induced by $(u,F)$. Since ${\rm Deg}(u)={\rm Deg}(g)$ and ${\rm Deg}(F)\le 4{\rm Deg}(g)$, \cite[Proposition IV.11.9]{FarkasKra1992} ensures that the total algebraic degree of $P$ satisfies ${\rm Deg} (P)\le 5{\rm Deg}(g)$. 
Setting $F_1=\Re F$ and $F_2=\Im F$, \eqref{eq:C-m_Omega} and \eqref{eq:C-m_Omega1} give
\begin{equation}\label{eq:Milnor-Thom}
C=\{(u,\imath F_2)\in\r\times\imath\r\colon P(u,\imath F_2)=0\}.
\end{equation} 
Write $Q_1(u,F_2)=\Re P(u,\imath F_2)$ and  $Q_2(u,F_2)=\Im P(u,\imath F_2)$. Note that $Q_1$ and $Q_2$ are real polynomials in two real variables with the total algebraic degree ${\rm Deg}(Q_j)\le {\rm Deg} (P)\le 5{\rm Deg}(g)$, $j=1,2$. Thus, \cite[Theorem 8.1]{Wallach1996} ensures that the algebraic system $Q_1(u,F_2)=Q_2(u,F_2)=0$  has at most $25{\rm Deg}(g)^2$ solutions in $\r^2$. 
Since ${\rm Deg}(u)={\rm Deg}(g)$, this and \eqref{eq:Milnor-Thom} imply  \eqref{eq:cuspglobal-FTC-C}, which proves \eqref{eq:cuspglobal-FTC}.

Finally, \eqref{eq:mOmega} and \eqref{eq:cuspglobal-FTC} give $\#(X^{-1}(L))\le 6{\rm Deg}(N)+\chi(M)+ 25{\rm Deg}(N)^3$.
\qed
%\smallskip

\noindent
{\em Proof of Theorem 1.3.} 
Follow \cite[proof of Theorem 1.3]{AlarconCastro-InfantesLopez2019CVPDE} but choosing each set $C_j$ $(j=1,2)$ consisting of $m+6r+25r^3+1$ points and applying Theorem 1.4.
\qed
%\smallskip

Note that \cite[Corollary 4.5]{AlarconCastro-InfantesLopez2019CVPDE} remains to hold true; just replace \cite[Theorem 1.4]{AlarconCastro-InfantesLopez2019CVPDE} by Theorem 1.4 in its proof. Finally, the statement of \cite[Corollary 4.4]{AlarconCastro-InfantesLopez2019CVPDE} has to be corrected by adding $50r^3$ points to the set $A^*_{r;m}$; here is its correct formulation.

%\smallskip
\noindent{\bf Corollary 4.4.}
{\em Let $r$ and $m$ be as in Theorem 1.3. There is a set $A^*_{r;m}\subset\r^3$, consisting of $12r+50r^3 +2m+2$ points,  such that if $X\colon M\to\r^3$ is a complete orientable immersed minimal surface with   $\chi(M)\le m$ and  $A^*_{r;m}\subset X(M)$, then the absolute value of the total curvature $|{\rm TC}(X)|>4\pi r$.}

\noindent{\em Acknowledgements.}
The authors were partially supported by the MINECO/FEDER grant no. MTM2017-89677-P, Spain.  I.\ Castro-Infantes was also partially supported by the MICINN/FEDER project PGC2018-097046-B-I00, and Fundaci\'on S\'eneca project 19901/GERM/15, Spain, and by the MICINN grant FJC2018-035533-I co-financed by the ESF.

%%%%%%%%%%
%%%%%%%%%%
%%%%%%%%%%
%%%%%%%%%%   THE BIBLIOGRAPHY
%%%%%%%%%%
%%%%%%%%%%

%{\bibliographystyle{abbrv} \bibliography{references}}

%%%%%%%%%%
%%%%%%%%%%
%%%%%%%%%%
%%%%%%%%%%   AFFILIATIONS
%%%%%%%%%%
%%%%%%%%%%

\bigskip

\noindent Antonio Alarc\'{o}n $\cdot$ Francisco J.\ L\'opez

\noindent Departamento de Geometr\'{\i}a y Topolog\'{\i}a e Instituto de Matem\'aticas (IEMath-GR), Universidad de Granada, Campus de Fuentenueva s/n, E--18071 Granada, Spain.

\noindent  e-mail: {\tt alarcon@ugr.es} $\cdot$ {\tt fjlopez@ugr.es}

\smallskip

\noindent Ildefonso Castro Infantes

\noindent Departamento de Matem\'aticas, Universidad de Murcia, Campus de Espinardo 30100, Murcia, Spain.

\noindent  e-mail: {\tt ildefonso.castro@um.es}

\end{document}